# 11 can be reduced to 10

Abraham Berman and Eliahu Levy

Technion – Israel Institute of Technology

**Abstract**

Laffey and Šmigoc proved that for every 2x2 doubly nonnegative integer matrix A, icpr(A) ≤ 11. We prove that 11 can be replaced by 10 and show that for many small matrices, even by 9.

Keywords: Complete positivity, complete positivity over the integers, sums of perfect squares, integer CP rank.

## Introduction

Laffey and Šmigoc (2019) proved that every 2x2 doubly nonnegative integer matrix can be decomposed as $BB^T$, where B is a 2x11 nonnegative integer matrix. They pointed out that 11 cannot be replaced by 8 and asked if it can be replaced by 9 or 10. We prove that 11 can be replaced by 10 and show that for many small matrices, even by 9.

The interest in the problem stems from the theory of completely positive matrices so we start the paper with a very short introduction to this theory. All the numbers (vectors and matrices) in the paper are nonnegative and after the section on complete positivity, they are nonnegative integers. The main tool used in the paper are the theorems of Lagrange and Legendre on expressing a number as sum of perfect squares, so we cite these important theorems. The proofs of the new results involve many special cases, so we divide them to description of the proof idea and tables that depict the special cases.

## Complete positivity

**Definition 1 (***Completely positive matrix and the CP rank***)**
A matrix A is *completely positive* if it can be decomposed as $A = BB^T$, where B is a nonnegative matrix. The smallest number of columns in such matrix B is called the *CP rank* of A.
 References on properties and applications of complete positivity and on bounds on the CP rank are [Berman and Shaked-Monderer, 2003] and [Shaked-Monderer and Berman, 2021].

**Definition 2 (***Doubly nonnegative matrix***)**
A positive semidefinite matrix that is also elementwise nonnegative is called *doubly nonnegative (DNN)*.

Obviously, completely positive matrices are doubly nonnegative, but the necessary condition is not sufficient.

**Definition 3** (*Matrix realization*)
Let G be a simple graph with vertices {1, 2, ..., n}. A symmetric matrix A is a matrix realization of G if for $i \neq j$, $a_{ij} \neq 0$ if and only if i and j are adjacent (there is no restriction on the diagonal entries of A).

**Definition 4** (*Completely positive graph*)
A simple graph G is *completely positive* if every doubly nonnegative matrix realization of G is completely positive.

**Theorem 1 (**Kogan and Berman,1993)
A graph G is completely positive if and only if it does not contain an odd cycle of length > 4.

**Remark 1**
An equivalent characterization is, [Trotter, 1979], that G is the line graph of a perfect graph.

**Corollary 1** (Maxfield and Minc, 1962)
Graphs with less than 5 vertices are completely positive.

## Integer completely positive matrices

**Definition 5** (*Completely positive matrix over the integers and the integer CP rank*)
A matrix A is *completely positive over the integers* if it can be decomposed as A = BB$^T$, where B is an integer nonnegative matrix.
The smallest number of columns in such matrix B is called the *integer CP rank* of A and is denoted by icpr(A).

The matrix $\begin{pmatrix} 1 & 1 & 1 \\ 1 & 2 & 0 \\ 1 & 0 & 2 \end{pmatrix}$ is an example of a 3x3 integer completely positive matrix that is not completely positive over the integers, but for n=2, Laffey and Šmigoc

(2018) proved that every 2x2 completely positive matrix (that by Corollary 1 is the same as every 2x2 DNN matrix) is completely positive over the integers. In 2019 they proved:

**Theorem 2**
For every 2x2 DNN matrix A, icpr(A) $\leq$ 11.

The matrix $\begin{pmatrix} 8 & 1 \\ 1 & 8 \end{pmatrix}$ shows that 11 cannot be replaced by 8. They asked if it can be replaced by 9 or 10. We prove that it can be replaced by 10 and show that for many small matrices, even by 9.

A key step in the proof of Theorem 2 was:

**Theorem 3**
If A = $\begin{pmatrix} a & b \\ b & c \end{pmatrix}$ is DNN and b>c (or b>a), then there exists a DNN matrix
A' = $\begin{pmatrix} a' & b' \\ b' & c' \end{pmatrix}$ , where $a' \geq b'$ and $c' \geq b'$, and A' and A have the same icpr.

## Sums of perfect squares

The central tool in the paper is the theorems of Lagrange and Legendre.

**Theorem 4** (Lagrange's theorem)
Every natural number can be written as sum of 4 perfect squares.

**Theorem 5** (Legendre's theorem)
A natural number can be written as sum of 3 perfect squares if and only if it is not of the form $4^k(8m+7)$.

**Corollary 2**
The numbers 8k+1, 8k+2, 8k+3, 8k+5, 8k+6 can be written as sum of 3 perfect squares.

**Definition 6** (*Good numbers and bad numbers*)
We say that a number x is *good* if x (mod 8) $\in$ {1, 2, 3, 5, 6} and *bad* if
x (mod 8) $\in$ {0, 4, 7}.

**Definition 7** (*bad, good and very good triplets*).

A triplet (e, f, g) is *good* if at least two of the numbers (e-f) (mod 8), f , (g-f) (mod 8) are good. It is *very good* if all three are good  and it is *bad* if at least two numbers are bad.

**Remark 2**

In terms of  integer CP rank, Theorem 4 says that for a 1x1 matrix, the icpr is bounded by 4 and for the numbers in Corollary 2 (the good numbers) it is less than or equal to 3.

## Conjectures and theorems

### Conjecture 1

Every 2x2 DNN matrix A can be decomposed as A = BB$^T$ where  B $\in$ R$^{2 \times 9}$.

To bound the CP rank, we observe that icpr(A+B) ≤ icpr(A)+ icpr(B).

Let A = $\begin{pmatrix} a & b \\ b & c \end{pmatrix}$ be DNN .

Without loss of generality, we can assume that c ≥ a. By Theorem 3, we can assume that c ≥a ≥b≥0.

If b=0 , icpr(A) = icpr (a) + icpr(c) ≤8, so we can assume that b>0 .

If a=b , then A = $\begin{pmatrix} 0 & 0 \\ 0 & c-b \end{pmatrix} + \begin{pmatrix} b & b \\ b & b \end{pmatrix}$ and icpr(A) ≤ 8 so we can assume that c  ≥ a > b > 0.

If b$^2$ ≤  c, then A = $\begin{pmatrix} a-1 & 0 \\ 0 & c-b^2 \end{pmatrix} + \begin{pmatrix} 1 & b \\ b & b^2 \end{pmatrix}$ and icpr(A) ≤ 9 so we can assume that   b$^2$ > c  ≥ a >  b > 0.

To prove the conjecture we suggest:

### Conjecture 2

For any two natural numbers a>b, there exist vectors u = (u$_1$ u$_2$ u$_3$ u$_4$ u$_5$) and v = (v$_1$ v$_2$ v$_3$ v$_4$ v$_5$)  such that uu$^T$ = a, and  vu$^T$ = b  and vv$^T$ ≤ b.

### Conjecture 3

For any two natural numbers a>b, there exist vectors u = (u$_1$u$_2$u$_3$u$_4$u$_5$) and v = (v$_1$v$_2$v$_3$v$_4$v$_5$) such that v ≤ u, uu$^T$ = a, and  vu$^T$ = b.

**Remark 3**

Obviously, Conjecture 3 implies Conjecture 2. Conjecture 2, in turn, implies Conjecture 1 with $\begin{pmatrix} u_1 & u_2 & u_3 & u_4 & u_5 & 0 & 0 & 0 & 0 \\ v_1 & v_2 & v_3 & v_4 & v_5 & w_1 & w_2 & w_3 & w_4 \end{pmatrix}$ where $c - \sum_{i=1}^{5} v_i^2 = \sum_{i=1}^{4} w_i^2$.

To prove the first new result, we need the following definitions and lemmas:

**Definition 8** (*Spanning vector*)

A vector $u = (u_k \ldots u_2\, u_1)$ is *spanning* if for every $b \leq uu^T$, there exists vector $v$ such that $v \leq u$ and $vu^T = b$.

**Definition 9** (*Step vector*)

A vector $u = (u_k \ldots u_2\, u_1)$ is a *step vector* if $u_1 = 1$ and $(u_{(i+1)} - u_i) \in \{0,1\}$, $i = 1,\ldots,k-1$.

**Lemma 1**

Step vectors are spanning.

This is a special case of:

**Lemma 2**

If $u = (u_k \ldots u_2\, u_1)$ is a spanning vector and $u_{(k+1)} \leq uu^T + 1$, then $(u_{(k+1)} \ldots u_2\, u_1)$ is also spanning.

**Proof**

Every number between 1 and $\sum_{i=1}^{k+1} u_i^2$ can be written as
$pu_{k+1} + q$, $0 \leq p \leq u_{k+1}$, $0 \leq q \leq \sum_{i=1}^{k} u_i^2$ ∎

**Theorem 6**

If $A = \begin{pmatrix} a & b \\ b & c \end{pmatrix}$ is DNN and $a \leq 64$, then $A$ can be decomposed as $A = BB^T$, where $B \in R^{2 \times 9}$.

**Proof**

We show that Conjecture 3 holds for $a \leq 64$.
In Table 1, we show that for every number between 1 and 64, except 33, there is a spanning vector. In Table 2, we construct vectors $u$ and $v$, for $a = 33$ and $b < 33$. ∎

**Remark 4**

65 does not have a spanning vector but for every b < 65 there are vectors u and v that satisfy Conjecture 3. Theorem 6 is used in the proof of Theorem 7 and 64 is sufficient for this purpose.

**Remark 5**

Conjecture 1 holds for every number that has a spanning vector. Unfortunately, the number of natural numbers that have a spanning vector is finite. The largest integer that has a spanning vector is $1^2 + 2^2 + 6^2 + 42^2 + 1806^2 = 3263441$.

**Theorem 7**

Every 2x2 DNN matrix A can be decomposed as $A = BB^T$ where $B \in R^{2 \times 10}$.

**Proof**

Let $A = \begin{pmatrix} a & b \\ b & c \end{pmatrix}$, $b^2 > c \geq a > b$.

We want to show that $icpr(A) \leq 10$.

If $a < 65$, then, by Theorem 6, $icpr(A) \leq 9$.

If $a > 64$ then $b > 8$ and $c > 64$.

Consider the remainders $\alpha, \beta, \gamma$ when a, b and c are divided by 8.

$$A = \begin{pmatrix} 8t + \alpha & 8r + \beta \\ 8r + \beta & 8s + \gamma \end{pmatrix} = (8r + \beta) J + \begin{pmatrix} 8(t-r) + \alpha - \beta & 0 \\ 0 & 8(s-r) + \gamma - \beta \end{pmatrix},$$

where $J = \begin{pmatrix} 1 & 1 \\ 1 & 1 \end{pmatrix}$.

If $(\alpha, \beta, \gamma)$ is very good then, by Corollary 2, $icpr(A) \leq 9$. If it is good then, by the same corollary, $icpr(A) \leq 10$. If $(\alpha, \beta, \gamma)$ is bad we subtract from A a rank one matrix $\begin{pmatrix} x \\ y \end{pmatrix} (x \quad y)$, where x, y, xy are less than or equal to 7. Since $a > 64$, $\begin{pmatrix} a & b \\ b & c \end{pmatrix} - \begin{pmatrix} x^2 & xy \\ xy & y^2 \end{pmatrix} \geq 0$.

Denoting $\alpha' = \alpha - x^2 \pmod{8}$, $\gamma' = \gamma - y^2 \pmod{8}$, $\beta' = \beta - xy \pmod{8}$, we get

$$A = \begin{pmatrix} x \\ y \end{pmatrix} (x \quad y) + (8r + \beta - xy) J + \begin{pmatrix} 8(t-r) + \alpha - \beta - x^2 & 0 \\ 0 & 8(s-r) + \gamma - \beta - y^2 \end{pmatrix} =$$

$$= \begin{pmatrix} x \\ y \end{pmatrix} (x \quad y) + (8r + \beta') J + \begin{pmatrix} 8(t-r) + \alpha' - \beta' & 0 \\ 0 & 8(s-r) + \gamma' - \beta' \end{pmatrix}$$

To complete the proof, we show in Table 3 that for every bad triplet (α, β , γ ) there is a pair (x, y) such that (α', β', γ') is very good. The pairs used are (1,2) (1,3) (1,4) (1,5) (1,6) (1,7) (2,1) (2,3) so $1 \leq x, xy, y \leq 7$. ∎

## Tables

### Table 1

For every number up to 64, except 33, there is a spanning vector (SV) that can be obtained by expressing n-1 as sum of 4 perfect squares.

| a | SV | a | SV | a | SV | a | SV |
|---|---|---|---|---|---|---|---|
| 64 | 73211 | 48 | 53321 | 32 | 52111 | 16 | 32111 |
| 63 | 7321 | 47 | 54211 | 31 | 5211 | 15 | 3211 |
| 62 | 72221 | 46 | 44321 | 30 | 4321 | 14 | 321 |
| 61 | 73111 | 45 | 6221 | 29 | 51111 | 13 | 31111 |
| 60 | 7311 | 44 | 54111 | 28 | 43111 | 12 | 3111 |
| 59 | 72211 | 43 | 62111 | 27 | 4311 | 11 | 22111 |
| 58 | 7221 | 42 | 6211 | 26 | 42211 | 10 | 2211 |
| 57 | 6421 | 41 | 44221 | 25 | 4221 | 9 | 221 |
| 56 | 63311 | 40 | 53211 | 24 | 33211 | 8 | 21111 |
| 55 | 54321 | 39 | 5321 | 23 | 3321 | 7 | 2111 |
| 54 | 63221 | 38 | 52221 | 22 | 4211 | 6 | 211 |
| 53 | 44421 | 37 | 53111 | 21 | 421 | 5 | 21 |
| 52 | 54311 | 36 | 5311 | 20 | 41111 | 4 | 1111 |
| 51 | 63211 | 35 | 52211 | 19 | 4111 | 3 | 111 |
| 50 | 54221 | 34 | 43221 | 18 | 3221 | 2 | 11 |
| 49 | 62221 | 33 | ****** | 17 | 22221 | 1 | 1 |

**Table 2**

For a=33 there is no spanning vector. We must use all three representations of 33 as sums of up to 5 squares, $33 = 4^2 + 4^2 + 1 = 5^2+2^2+2^2 = 4^2+3^2+2^2+2^2$.

In the table, we denote u=(44100) by 441, u=(52200) by 522 and u= (43220) by 4322 and construct (shortened) vector v for $1 \leq b \leq 32$.

| b | 441 | b | 522 | b | 4322 |
|---|---|---|---|---|---|
| 1 | 001 | 2 | 001 | 3 | 0100 |
| 4 | 010 | 6 | 021 | 18 | 4001 |
| 5 | 011 | 7 | 101 | 30 | 4222 |
| 8 | 020 | 10 | 200 | | |
| 9 | 021 | 11 | 121 | | |
| 12 | 030 | 14 | 202 | | |
| 13 | 031 | 15 | 300 | | |
| 16 | 040 | 19 | 311 | | |
| 17 | 041 | 22 | 401 | | |
| 20 | 140 | 23 | 322 | | |
| 21 | 141 | 26 | 421 | | |
| 24 | 240 | 27 | 501 | | |
| 25 | 241 | 31 | 521 | | |
| 28 | 340 | | | | |
| 29 | 341 | | | | |
| 32 | 440 | | | | |

**Table 3**

We consider the bad triplets (α, β, γ). The symmetry allows us to consider only the cases γ ≤ α. We list them according to the pairs (x, y) used in the proof of Theorem 7. The columns are organized in decreasing order, first of β, then α and then γ.

**x = 1, y = 2**

α  7 7 7 7 7 4 3 3 5 5 5 1 6 6 6 5 5 5 4 4 4 4 1 0 7 7 7 3 3 4 4 2 1 0
β  7 7 7 7 7 7 7 7 5 5 5 5 4 4 4 4 4 4 4 4 4 4 4 4 3 3 3 3 3 0 0 0 0 0
γ  7 6 4 3 2 3 3 2 5 4 1 1 4 3 0 4 3 0 4 3 1 0 0 0 7 3 2 3 2 4 3 0 0 0

**x = 1, y = 3**

α  6 6 6 5 5 2 4  4 7 7 7 3 3  5 5 5 1 1 0 7 7 7 7 4 4 3
β  6 6 6 6 6 6 5  5 4 4 4 4 4  1 1 1 1 1 1 0 0 0 0 0 0 0
γ  6 5 2 5 2 2 4  1 4 3 0 3 0  5 1 0 1 0 0 7 3 1 0 1 0 0

**x = 1, y = 4**

α  7 7 7 6 6 6 6 6 2 1
β  7 7 7 7 7 7 7 7 2 2
γ  5 1 0 6 5 4 1 0 1 1

**x = 1, y = 5**

α  3 3 2 2 7 7 7
β  4 4 4 3 0 0 0
γ  2 1 3 6 5 2

**x = 1, y = 6**

α  5 3 6  6 5 5
β  7 7 0  0 0 0
γ  3 0 4  0 4 0

**x = 1, y = 7**

α   6 6 6 6 6 2 7
β  7 7 2 2 2 2 0
γ  3 2 6 2 1 2 4

**x = 2, y = 1**

α  3 4
β  7 0
γ  1 2

**x = 2, y = 3**

α    4
β    4
γ    2